\documentclass[letterpaper, 10 pt, conference]{ieeeconf}  % Comment this line out
\usepackage{epsfig} % for postscript graphics files
\usepackage{mathptmx} % assumes new font selection scheme installed
\usepackage{times} % assumes new font selection scheme installed
\usepackage{amsmath} % assumes amsmath package installed
\usepackage{amssymb}  % assumes amsmath package installed
\usepackage{graphics} % for pdf, bitmapped graphics files
\usepackage{wrapfig}	% enables floating wrapped figures

\usepackage[]{algorithm2e}
\usepackage{float}

\IEEEoverridecommandlockouts                              % This command is only
\title{\LARGE \bf
Stochastic Nonlinear Model Predictive Control with State Estimation by Incorporation of the Unscented Kalman Filter
}

\author{Eric Bradford\textsuperscript{1} and Lars Imsland\textsuperscript{2}% <-this % stops a space
\thanks{1 E. Bradford is with the Faculty of Information Technology and Electrical Engineering, Engineering Cybernetics,
        NTNU, 7491 Trondheim, Norway
        {\tt\small eric.bradford@ntnu.no}}%
\thanks{2 L. Imsland is with the Faculty of Information Technology and Electrical Engineering, Engineering Cybernetics,
        NTNU, 7491 Trondheim, Norway
        {\tt\small lars.imsland@ntnu.no}}%
}

\begin{document}

\maketitle
\thispagestyle{empty}
\pagestyle{empty}

%%%%%%%%%%%%%%%%%%%%%%%%%%%%%%%%%%%%%%%%%%%%%%%%%%%%%%%%%%%%%%%%%%%%%%%%%%%%%%%%
\begin{abstract}
Nonlinear model predictive control has become a popular approach to deal with highly nonlinear and unsteady state systems, the performance of which can however deteriorate due to unaccounted uncertainties. Model predictive control is commonly used with states from a state estimator in place of the exact states without consideration of the error. In this paper an approach is proposed by incorporating the unscented Kalman filter into the NMPC problem, which propagates uncertainty introduced from both the state estimate and additive noise from disturbances forward in time. The feasibility is maintained through probabilistic constraints based on the Gaussian approximations of the state distributions. The concept of "robust horizon" is introduced to limit the open-loop covariances, which otherwise grow too large and lead to conservativeness and infeasibility of the MPC problem. The effectiveness of the approach was tested on a challenging semi-batch reactor case study with an economic objective.   
\end{abstract}

%%%%%%%%%%%%%%%%%%%%%%%%%%%%%%%%%%%%%%%%%%%%%%%%%%%%%%%%%%%%%%%%%%%%%%%%%%%%%%%%
\section{INTRODUCTION}
Model predictive control (MPC) was developed in the late seventies as a method to deal with system constraints and strongly coupled, multivariable plants. MPC is the only advanced control approach that has been applied in industry in a large fashion \cite{Maciejowski2002}. MPC solves at each sampling instance an open-loop, optimal control problem (OCP) based on an explicit process model to determine a finite sequence of control actions to take. The first of these control actions is implemented, while discarding the rest \cite{Mayne2000}. Feedback is implicitly introduced in this process by the state and bias update using the measurements available at each sampling time \cite{Mesbah2016}. 

MPC methods based on linear models have found a multitude of successful applications in industry, in particular in the process industry \cite{Camacho2012}. Linear MPC (LMPC) theory is relatively mature and well-established in practice. LMPC however is inadequate to handle processes with strong nonlinearities or at unsteady state, such as batch processes. In addition, higher productivity demands and tighter environmental regulations require a more accurate description of the plant and hence motivate the use of nonlinear models \cite{Findeisen2003}. In addition, nonlinear MPC (NMPC) allows optimization with respect to economic criteria directly \cite{Rawlings2009}. 

The introduction of uncertainty in the system can lead to sub-optimal behaviour and failures of the MPC algorithm. Key questions are the maintenance of stability, constraint satisfaction and recursive feasibility for the uncertainty in question. The solution of MPC is often close to its constraints and hence can easily lead to infeasibilities due to unaccounted uncertainties. Therefore, the development of MPC approaches that make use of an explicit description of the uncertainties has been of major interest over the past two decades \cite{Mayne2000}.

Robust NMPC (RNMPC) describe a series of methods that assume uncertainties to be given by bounded sets. Prominent techniques include min-max \cite{Campo1987}, tube-based \cite{Mayne2011} and multi-stage NMPC \cite{Lucia2014}. Min-max NMPC formulations focus on minimizing cost while satisfying the constraints under the worst-case realization. Open-loop min-max methods have been shown to be inadequate to deal with the spread of state trajectories, while closed-loop min-max approaches are difficult to solve \cite{Scokaert1998}. Tube-based NMPC was subsequently devised to address the limitations of min-max techniques. Tube-based NMPC aims to keep the trajectories in a tube that is computed offline. The tube is centred around a nominal trajectory, while the so-called "ancillary" controller determines a control policy that ensures that the trajectory of the real uncertain system remains in the tube \cite{Mayne2011}. Lastly, a multi-stage NMPC approach has been suggested in which the uncertainty is modelled by a scenario tree approach from stochastic programming. The method can explicitly take into account information available through feedback at each sampling time. However, the procedure quickly becomes intractable, since the size of the optimization problem scales exponentially with the time horizon, number of uncertainties and uncertainty levels \cite{Lucia2014}. 

An alternative to RNMPC is stochastic NMPC (SNMPC) in which uncertainties are given by known probability distributions. SNMPC allows to systematically trade-off the conservativeness of a solution by addressing constraints probabilistically and allowing an admissible level of constraint violation. The majority of work for stochastic MPC has been carried out for linear MPC, while SNMPC has not received much attention \cite{Mesbah2016}. A simple solution to SNMPC can be found by \cite{Cannon2009a}, who linearises the nonlinear system successively and then applies a probabilistic tube method. A popular approach in SNMPC is given by the use of polynomial chaos (PC) expansions, which is a computationally efficient tool for accelerating sampling-based techniques \cite{Streif2014}. In this method, implicit mappings between variables/parameters and the states are replaced by orthogonal polynomials. A disadvantage of this approach is that the computational cost scales exponentially with the number of uncertainty parameters. Apart from PC, importance sampling methods have been put forward \cite{Maciejowski2007}, which are generally more efficient than standard Monte-Carlo methods, but do not take gradient information into account. \cite{Goodwin2009} uses the multistage stochastic programming approach with a novel scenario generation method to solve the SNMPC problem taking imperfect feedback information into account, similar to \cite{Lucia2014}. The approach does, however, quickly become intractable due to the aforementioned scaling.  

In this paper a method is proposed that yields a tractable approximation to the SNMPC problem based on previous work by \cite{Yan2005}, who incorporated the Kalman filter into LMPC. Further, in \cite{Liu2014} and \cite{Volz2016} the unscented transformation is used to propagate additive disturbance error for a nonlinear system, which is shown to perform well. In \cite{Liu2014} proofs are given for recursive feasibility. We propose to incorporate the unscented Kalman filter (UKF) into the MPC problem to estimate and propagate the mean and covariance of the states along the time horizon, which takes into account both error due to additive noise from disturbances and measurement error and also error introduced due to the imperfect knowledge of the states through the state estimation, which was not covered in \cite{Liu2014} or \cite{Volz2016}. Further in \cite{Liu2014} and \cite{Volz2016} the problem of growing covariances was not addressed. The problem of growing covariances in this paper is addressed by introduction of the "robust-horizon", which is a cheap solution to the issue. The closest to the here proposed method can be found in \cite{Farrokhsiar2012}, in which the UKF is also incorporated online for a specific control problem of mobile robots. In \cite{Farrokhsiar2012} the covariances are kept limited by predicting future measurements, which however is expensive compared to the solution given here of the "robust horizon". The previous work has focused solely on the application of the UKF in MPC for regulatory problems. We test it on a challenging economic MPC problem based on a semi-batch reactor. 

The paper is structured as follows. In the second section a general formulation is given of the SNMPC problem we wish to solve. In the next section the UKF is introduced and incorporated into the MPC problem to yield the algorithm. Further, a simple solution is given to the problem of increasing covariances. In the fourth section a challenging case study is introduced to control a semi-batch reactor with an economic objective. In the fifth section the results of the case study are illustrated and discussed. Finally, in the last section we draw some conclusions and propose future work to be carried out.   

\section{NONLINEAR MODEL PREDICTIVE CONTROL WITH CHANCE CONSTRAINTS}
In this report we consider a general discrete-time stochastic nonlinear system with additive noise, described by:
\begin{align}
x(k+1)&=f(x(k),u(k)) + w(k) \label{eq:eq1} \\
y(k)&=h(x(k),u(k)) + \nu(k) \label{eq:eq2}
\end{align}
where $x \in \mathbb{R}^{n_x}$ denotes the system states, $u \in \mathbb{R}^{n_u}$ represents the control inputs and $y \in \mathbb{R}^{n_y}$ denotes the system measurements; the additive disturbance term $w$ lies in $\mathbb{R}^{n_x}$ and the additive measurement noise $\nu$ lies in $\mathbb{R}^{n_y}$. The equations $f:\mathbb{R}^{n_x} \times \mathbb{R}^{n_u} \rightarrow \mathbb{R}^{n_x}$ and $h:\mathbb{R}^{n_x} \times \mathbb{R}^{n_u} \rightarrow \mathbb{R}^{n_y}$ represent the system dynamics of the states and the measurements respectively. 

It is assumed that $\{w(k)\}$ and $\{v(k)\}$ are sequences of zero mean normal independent random variables with variances $\Sigma_w(k)$ and $\Sigma_v(k)$ at stage $k$ respectively. Furthermore, the prior density of $x(0)$ is assumed to be normal with mean $\hat{x}_0$ and variance $\Sigma_{x_{0}}$. Let $\mathcal{Y}_n$ represent the sequence of measurements collected up to time $n$, $\mathcal{Y}_n = \{(u(i),y(i))\}_{i=1,\ldots,n}$. Then, the notations $\mathbb{E}_{\mathcal{Y}_n}(\cdot)$ and $\mathbb{P}_{\mathcal{Y}_n}(\cdot)$ refer to the conditional expectation and probability respectively conditioned on $\mathcal{Y}_n$ \cite{Yan2005}. The aim of the SNMPC algorithm is to adjust the probability distributions of the future states in the time horizon to lie within predefined constraints and give an optimum performance with respect to the objective with imperfect information available through ${\mathcal{Y}_n}$ at stage $n$. The general chance constrained, finite-horizon SNMPC problem at time $n$ we consider is given by:
\\ \vspace{-5pt}

\textbf{Finite-horizon SNMPC problem with chance constraints}
\begin{equation}
\begin{aligned}
& \underset{\mathbf{u}_N}{\text{minimize}} \quad \mathbb{E}_{\mathcal{Y}_n}(J(N,x(n),\mathbf{u}_N)) \\
& \text{subject to}  \\
& x(n+k+1)=f\left(x(n+k),u(n+k)\right)+w(n+k) \\
& y(n+k)=h(x(n+k),u(n+k))+\nu(n+k) \\
& \mathbb{P}_{\mathcal{Y}_n}(x(n+k) \in \mathbb{X}_k) \geqslant p_k \quad \forall k \in \{1,...,N\} \\
& u(n+k) \in \mathbb{U}_k \quad \forall k \in \{0,...,N-1\} \\
\end{aligned}
\label{eq:eq3}
\end{equation}
where $N$ represents the length of the time horizon, $\mathbf{u}_N:=\{u(n),\ldots,u(n+N-1)\}$ denotes the decision variables over the finite-horizon $N$ from an initial stage $n$, $\mathbb{X}_k \subset \mathbb{R}^{n_x}$ denotes a compact set of state constraints, $\mathbb{U}_k \subset \mathbb{R}^{n_u}$ are the compact sets of input constraints and $J(N,x(n),\mathbf{u}_N))$ the probabilistic objective function. The joint probability constraints can be violated by only a specific rate, given by $p_k \in (0,1) \subset \mathbb{R}$.       

The finite-horizon OCP given in (\ref{eq:eq3}) is based on both the dynamics of the states given in (\ref{eq:eq1}) and the dynamics of the measurements given in (\ref{eq:eq2}). The problem considers joint probability constraints on the states and deterministic inequality constraints on the inputs.   

\section{INCORPORATION OF THE UNSCENTED KALMAN FILTER}
\subsection{Unscented Kalman filter with additive noise}
The complexity of the constrained stochastic optimization problem in (\ref{eq:eq3}) is prohibitive, since it would require the full determination and propagation of the entire conditional distribution of the states through nonlinear transformations \cite{Yan2005}. To make progress we therefore need to make assumptions. It is assumed that the states follow a Gaussian distribution and hence we only need to predict and propagate the mean and the covariance of the states. 

Let $\hat{x}(n+k|n)$ be the mean and $\Sigma_x(n+k|n)$ be the variance of the state vector at time $n+k$ given data $\mathcal{Y}_n$. We require a method to find $\hat{x}(n+k|n)$ and $\Sigma_x(n+k|n)$ from an initial time $n$ up to time $k=N$. Bayesian recursive filtering deals with the problem of propagating probability distributions given a set of observations. For linear systems the finite-dimensional Kalman filter can be used to directly propagate the Gaussian distributions, while for nonlinear systems the probability distributions need to be approximated at each stage. In this report we use the UKF to approximate the probability distributions of the states at each stage. The UKF for the case of additive noise for (\ref{eq:eq1}) and (\ref{eq:eq2}) to approximate $\hat{x}(n+k|n)$ and $\Sigma_x(n+k|n)$ can be stated as follows \cite{Wan2000}: \\ \vspace{-5pt}

\textbf{UKF with additive noise} \\
\vspace{4pt}
\text{Definition of Sigma-points}
\begin{multline}
\mathcal{X}(n+k-1|n) = \lbrack \hat{x}(n+k-1|n) \\ \hat{x}(n+k-1|n) +\sqrt{L+\lambda} \, \Sigma^{1/2}_x(n+k-1|n) \\ \hat{x}(n+k-1|n)-\sqrt{L+\lambda}  \, \Sigma^{1/2}_x(n+k-1|n) \rbrack \qquad
\label{eq:eq4}
\end{multline}
\vspace{4pt}
\text{Covariance and mean approximation of predictions}
\begin{subequations}
\begin{flalign}
& \mathcal{X}^{(i)}(n+k|n) = f(\mathcal{X}^{(i)}(n+k-1|n),u(n+k-1)) \\
& \hat{x}(n+k|n)=\sum_{i=0}^{2L}\omega^{\mu}_i \mathcal{X}^{(i)}(n+k|n) \\
\begin{split}
& \Sigma_x(n+k|n) = \sum_{i=0}^{2L}\omega^{c}_i(\mathcal{X}^{(i)}(n+k|n)-  \\ & \hat{x}(n+k|n))(\mathcal{X}^{(i)}(n+k|n)-\hat{x}(n+k|n))^T + \Sigma_w(n+k)
\end{split}
\end{flalign}
\label{eq:eq5}
\end{subequations}
\vspace{4pt}
\text{Covariance and mean approximation of observations}
\begin{subequations}
\begin{align}
& \phi^{(i)}(n|n-1) = h(\mathcal{X}^{(i)}(n|n-1),u(n-1)) \\
& \hat{y}(n|n-1) = \sum_{i=0}^{2L}\omega^{\mu}_i \phi^{(i)}(n|n-1) \\
\begin{split}
& \Sigma_{yy}(n|n-1) = \Sigma_v(n-1) + \sum_{i=0}^{2L}\omega^{c}_i (\phi^{(i)}(n|n-1)- \\ & \hat{y}(n|n-1))(\phi^{(i)}(n|n-1)-\hat{y}(n|n-1))^T
\end{split} \\
\begin{split}
& \Sigma_{xy}(n|n-1) = \sum_{i=0}^{2L}\omega^{c}_i (\mathcal{X}^{(i)}(n|n-1)- \\ & \hat{x}(n|n-1))(\phi^{(i)}(n|n-1)-\hat{y}(n|n-1))^T
\end{split}
\end{align}
\label{eq:eq6}
\end{subequations}
\vspace{4pt}
\text{Update of states from available measurements}
\begin{subequations}
\begin{align}
& K(n) = \Sigma_{xy}(n|n-1)\Sigma_{yy}(n|n-1)^{-1} \\
& \hat{x}(n|n) = \hat{x}(n|n-1) + K(n)(y(n)-\hat{y}(n|n-1)) \\
& \Sigma_x(n|n) = \Sigma_x(n|n-1) - K(n)\Sigma_{yy}(n|n-1)K(n)^T
\end{align}
\label{eq:eq7}
\end{subequations}
where $\mathcal{X}$ is a matrix of Sigma points with $\mathcal{X}^{(i)}$ corresponding to the rows of the matrix representing the Sigma points, $\phi^{(i)}$ are similarly the rows of transformed Sigma points of matrix $\phi$ to estimate the measurement mean and covariance, $\Sigma_{yy}$ denotes the covariance matrix of $y$, $\Sigma_{xy}$ represents the covariance matrix between $x$ and $y$, $\hat{y}$ the mean of $y$, $K$ the Kalman filter gain. Lastly, there are parameters $L$, $\lambda$, $\omega_i^m$ and $\omega_i^c$ that can be set as follows \cite{Wan2000}: \\
\begin{subequations}
\begin{align}
& L = 2n_x + 1 \\
& \lambda = \alpha^2(L+\kappa)-L \\
& \omega_0^m = \frac{\lambda}{L+\lambda} \\
& \omega_0^c = \frac{\lambda}{L+\lambda} + (1 - \alpha^2+\beta) \\
& \omega_i^m = \omega_i^c = \frac{1}{2(L+\lambda)} \quad \forall i  \in \{1,\ldots,2L-1\}
\end{align}
\label{eq:eq8}
\end{subequations}
where common values of $\alpha$, $\beta$ and $\kappa$ are $1e-3$, $2$ and $0$ respectively. 

\subsection{Robust horizon}
From (\ref{eq:eq5}c) we can see that the predicted conditional covariance usually increases with $k$ and hence the method becomes increasingly conservative with respect to the length of the time horizon. This can be seen by noting that the first part of the equation propagates the previous covariance forward, while the second part adds $\Sigma_w(n+k)$ on top each time. As pointed out in \cite{Yan2005} it is a conflict that larger time horizons lead to more and more difficult to solve MPC problems, while commonly large time horizons are associated with improved dynamic properties. Eventually for longer time horizons the MPC problem will become infeasible. In \cite{Hokayem2012} the control actions are replaced by parametrised feedback control laws to overcome this problem, which is computationally expensive however. To solve this problem we instead suggest introducing a so-called "robust horizon", similar to \cite{Lucia2014}, up to which the covariance matrix is propagated according to (\ref{eq:eq5}c) and after which the covariance matrix is kept the same. The rationale behind this step comes from the fact that the actual MPC controller implemented online has reduced covariances by the state and bias update through the measurements available at each sampling instance, which is otherwise not considered by the open-loop formulation. Hence, the following equation is added to the MPC problem:
\begin{equation}
\Sigma_x(n+k|n) = \Sigma_x(n+k-1|n) \quad \forall k\in \{t_R+1,\ldots,N\}
\label{eq:eq9}
\end{equation}
where $t_R$ is the time length of the robust horizon. 

\subsection{Formulation of UKF-SNMPC}
Dealing with general nonlinear, joint chance constraints with respect to the states is difficult. However, due to the fact that the underlying distributions of the states are Gaussian it is possible to give explicit expressions for linear joint probability constraints of the states. Therefore, in the proposed UKF-SNMPC algorithm the general probability constraints in (\ref{eq:eq3}) are replaced by probability constraints in the following form:
\begin{equation}
\mathbb{P}_{\mathcal{Y}_n}(H_k^Tx(n+k|n) \leq g_k) \geq p_k 
\label{eq:eq10}
\end{equation}
where $H_k \in \mathbb{R}^{n_x \times n_g}$ is a matrix representing the linear constraints and $g_k \in \mathbb{R}^{n_g}$ is a vector denoting the corresponding upper bounds. 

It was shown by \cite{VanHessem2002} using inscribed conic sets that the constraint in (\ref{eq:eq10}) can be given by the following constraints:
\begin{multline}
\Phi^{-1}(p_k) \sqrt{h_k^{(j)}{^T} \Sigma_x(n+k|n) h_k^{(j)}} + h_k^{(j)}{^T} \hat{x}(n+k|n) \leq g_k^{(j)} \\ \quad \forall j \in \{1,\ldots,n_{g_k}\} 
\label{eq:eq11}
\end{multline}
where $\Phi^{-1}(\cdot)$ is the quantile function of the standard Gaussian probability distribution, $n_{g_k}$ the number of rows of $g_k$, $h_k^{(j)}$ the $j^{th}$ row of $H_k$ and $g_k^{(j)}$ the corresponding $j^{th}$ value of $g$.  

Considering linear constraints of the form in (\ref{eq:eq10}), we can state the OCP for the SNMPC with incorporated UKF as follows: \\

\vspace{-5pt}
\textbf{Finite-horizon SNMPC problem with incorporated UKF and chance constraints}
\begin{equation}
\begin{aligned}
& \underset{\mathbf{u}_N}{\text{minimize}} \quad \mathbb{E}_{\mathcal{Y}_n}(J(N,x(n),\mathbf{u}_N)) \\
& \text{subject to} \\
& \mathbb{P}_{\mathcal{Y}_n}(H_k^Tx(n+k|n) \leq g_k) \geq p_k \quad \forall k \in \{1,...,N\} \\
& u(n+k) \in \mathbb{U}_k \quad \forall k \in \{0,...,N-1\} \\
& (\ref{eq:eq4}),(\ref{eq:eq5}),(\ref{eq:eq6}),(\ref{eq:eq7}), (\ref{eq:eq8}), (\ref{eq:eq9}) \\
\end{aligned}
\label{eq:eq12}
\end{equation}
where the probability constraints can be reformulated as shown in the previous section.

We can see that the overall algorithm includes both equations for state estimation and uncertainty propagation. (\ref{eq:eq5}b) is repeatedly used to propagate the mean of the states forward, while (\ref{eq:eq5}c) propagates the covariances up to a fixed "robust horizon", $t_R$, before the covariances are fixed according to (\ref{eq:eq9}). (\ref{eq:eq6}) and (\ref{eq:eq7}) are used only once in each open-loop problem to estimate the state according to the measurements available at time $n$ and to update the prior mean and covariance matrix. Therefore, to initialize the algorithm the previous mean, covariance matrix, control action and current measurement need to be supplied to the algorithm, the same as would be required for state estimation. The assumed error of the state estimate is given by the updated covariance matrix, which is propagated forward. The problem in (\ref{eq:eq12}) can be implemented in a receding-horizon fashion to yield a SNMPC controller. 

\section{BATCH REACTOR CASE STUDY}
\subsection{Dynamic model equations}
The use of the UKF-SNMPC algorithm is illustrated on the operation of a semi-batch reactor, which was adopted from \cite{Fogler1999}. The algorithm is applied to an economic MPC formulation, which aims to maximize the desired product $C$. The following series reactions take place inside the reactor, catalyzed by $H_2SO_4$:
\begin{equation*}
2A \xrightarrow[(1)]{k_{1A}}B \xrightarrow[(2)]{k_{2B}} 3C 
\end{equation*} 
The reactions taking place are all first-order with respect to the reactant concentration, however reaction $(1)$ is exothermic, while reaction $(2)$ is endothermic. The reactor is fitted with a heat exchanger. The control variables are given by the flow rate of pure A entering the reactor and the temperature of the heat exchanger. The evolution of the concentrations of A, B and C can be described by the following nonlinear differential algebraic equation (DAE) system:
\begin{subequations}
\begin{flalign}
& \dot{C}_A = -k_{1A}C_A + (C_{A0}-C_A)\frac{F}{V}, \\
& \dot{C}_B = 0.5k_{1A}C_A - k_{2B}C_B - C_B\frac{F}{V}, \\
& \dot{C}_C = 3k_{2B}C_B - C_C\frac{F}{V}, \\
\begin{split}
& \dot{T} = \frac{(UA(T_a-T)-FC_{A0}C_{P_A}(T-T_0)}{(C_A C_{P_A} +  C_B C_{P_B} + C_C C_{P_C})V + N_{H_{2}SO_4}C_{P_{H_{2}SO_4}}} + \\
& \frac{(-\Delta H_{Rx1A}k_{1A}C_A-\Delta H_{Rx2B}k_{2B}C_B)V}{(C_A C_{P_A} +  C_B C_{P_B} + C_C C_{P_C})V + N_{H_{2}SO_4}C_{P_{H_{2}SO_4}} }, 
\end{split} \hspace{-0.25in}  \\
& \dot{V} = F \\
& k_{1A} = A_1 \exp \left(-E_{1A}\left(\frac{1}{320}-\frac{1}{T}\right)\right) \\
& k_{2B} = A_2 \exp \left(-E_{2B}\left(\frac{1}{300}-\frac{1}{T}\right)\right)
\end{flalign}
\label{eq:eq13}
\end{subequations}
where $C_A$, $C_B$, $C_C$ are the concentrations in $\text{moldm}^{-3}$ of species $A$, $B$ and $C$ respectively, $T$ is the temperature in $K$ of the reactor and $V$ is the liquid volume in $dm^{-3}$. The parameters were kept at their nominal values, which can be found in \cite{Fogler1999}.

In compact form we can write $x=[C_A,C_B,C_C,T,V]^T$ and $u=[F,T_a]^T$. The continuous-time DAE system in (\ref{eq:eq13}) can be transformed to discrete-time using any numerical discretization, such as the Euler method. The discrete-time equation system can then be given by:
\begin{equation}
 x(k+1)=f(x(k),u(k)) + w(k) 
\label{eq:eq14}
\end{equation}
where $f(x(k),u(k))$ describes the DAE system in (\ref{eq:eq13}) and $w(k)$ is additive Gaussian noise with a constant covariance matrix $\Sigma_w = \text{diag}(1e-4,1e-4,2e-4,1,2)$.  

Lastly, the measurement dynamics need to be defined, which are given by the following simple equation:
\begin{equation}
y(k) = \begin{bmatrix}
    1 & 0 & 0 & 0 & 0 \\
    0 & 1 & 0 & 0 & 0 \\
    0 & 0 & 0 & 0 & 1 \\
  \end{bmatrix} x(k) + \nu (k)
\label{eq:eq15}
\end{equation}
where $\nu (k)$ is additive Gaussian noise with a constant covariance matrix given by $\Sigma_{\nu} =\text{diag}(1e-3,1e-3,1e-2)$. 

The measurement equation tells us that the variables $A$, $B$ and $V$ can be directly measured with additive noise, while measurements of $C$ and $T$ are not available.

\subsection{SNMPC problem}
In this section an OCP problem is defined based on an economic objective, which is subsequently solved in a receding-horizon fashion to yield a valid SNMPC for the system. The objective of the OCP problem given in (\ref{eq:eq16}) is to maximize the expected amount of $C$ with a penalty term added for excessive control actions. The probability constraints cover both path and terminal constraints. The volume is constrained to lie below $750 \text{dm}^{-3}$ and the temperature of the reactor is set to a safety limit of $440K$ for the entire time horizon. A terminal constraint was set for the concentration of reactant $A$ to lie below $0.5\text{moldm}^{-3}$. For all constraints the probability of constraint violation was set to 0.1. The flow rate is constrained to lie between $0\text{dm}^{-3}h^{-1}$ and $250\text{dm}^{-3}h^{-1}$, while the temperature of the heat exchanger is set to lie between $200K$ and $500K$. This is given by the following OCP problem:  

\begin{equation}
\begin{aligned}
& \underset{\mathbf{u}_N}{\text{minimize}} \quad -\mathbb{E}_{\mathcal{Y}_n}(x_2(n+N|n)x_4(n+N|n)) + \Delta U^T S \Delta U  \\
& \text{subject to}  \\
& \mathbb{P}_{\mathcal{Y}_n}(H_k^Tx(n+k|n) \leq g_k) \geq p_k \quad \forall k \in \{1,...,N\} \\
& u(n+k) \in [0,250]\times[200,500] \quad \forall k \in \{0,...,N-1\} \\
& (\ref{eq:eq4}),(\ref{eq:eq5}),(\ref{eq:eq6}),(\ref{eq:eq7}), (\ref{eq:eq8}), (\ref{eq:eq9}) \\
\end{aligned}
\label{eq:eq16}
\end{equation}
where $\Delta U = [u(n+k)-u(n+k-1)]_{k \in \{1,\ldots, N-1\}}$, $S = \text{diag}(2e-4,5e-5)$, $H_{1,\ldots,N-1} = \text{diag}(0,0,0,1,1)$, $g_{1,\ldots,N-1} = [0,0,0,440,750]$, $H_N = \text{diag}(1,0,0,1,1)$, $g_N = [0.5,0,0,440,750]$, $p_k=p_f=0.9$ and $\mathbb{E}_{\mathcal{Y}_n}(x_2(n+N|n)x_4(n+N|n))= \vspace{2pt} \hat{x}_2(n+N|n)\hat{x}_4(n+N|n)+\Sigma_x{_{2,4}}(n+N|n)$.

\subsection{Implementation}
The problem given in (\ref{eq:eq15}) is a standard OCP that is solved repeatedly at each new sampling time to yield an SNMPC controller. Each time it needs as input the previous control action, state estimate, covariance matrix and current measurement due to the incorporation of the state estimator. The OCP is solved using CasADi \cite{Andersson2013} by employing direct Collocation. The Collocation points were placed according to Radau quadrature rule and the degree of the polynomials was set to 3. For each control interval the state trajectories were approximated by two polynomials. The NLP  problem is solved utilising IPOPT \cite{Wachter2006}, which applies first and second order derivative information generated by automatic differentiation of CasADi to solve the NLP problem efficiently. The "real" nonlinear equation system was simulated using IDAS \cite{Hindmarsh2005}, which uses a backward differentiation formula implicit integration scheme. The random noise was generated pseudo-randomly from a Gaussian distribution. The computational work was carried out in Python on a Dell XPS 15 notebook with a Quad-core 6th Generation Intel i-7 processor with up to 3.5 GHZ and 16 GB of RAM. 

To show the effectiveness of the approach the case study was simulated 50 times, which leads each time to different trajectories given the uncertainty introduced through the disturbances and measurements. The parameters in (\ref{eq:eq4}-\ref{eq:eq9}) for the UKF  were set to the following values: $\alpha=0.4$, $\beta=2$, $\kappa = 0.1$. The length of the time horizon $N$ was set to $30$ with a time-interval length of $4/30\text{h}$. The complete simulation time was set to $6h$, which corresponds to 45 control inputs. At time $n=0$ the algorithm needs to be initiated by supplying it with the "previous" covariance matrix, mean and control action. These were set to $\hat{x}_0 = [0,0,0,290,100]$, $\Sigma_{x_0}=\text{diag}(1e-4,1e-4,1e-4,0.5,1)$ and $u_0=0$. The corresponding measurement that is required is obtained from ($\ref{eq:eq15}$). The variable $\hat{x}_0$ was perturbed by noise according to $\Sigma_{x_0}$ for each simulation. The robust horizon, $t_R$ was set to $2$.

\section{SIMULATION RESULTS}
In this section the semi-batch reactor simulation results are given. The trajectories for the states are shown in Fig. 1 to Fig. 5. Fig. 4 is of particular interest, since the temperature was a difficult constraint to adhere to given the steepness of the initial rise. As we can see however, the method is able to effectively reduce the constraint violation to near zero. In Fig. 5 the volume trajectories are shown, which are kept well below the constraint at $750dm^{-3}$. This indicates that the result is relatively conservative. The results could be potentially improved by using a robust horizon of length $1$ by allowing for more constraint violation. 

\vspace{-10pt}
\begin{figure}[H]
    \includegraphics[width=0.35\textwidth]{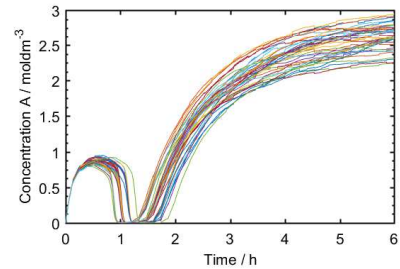}
    \centering
    \vspace{-10pt}
    \caption{Concentration A trajectories for 50 realizations}
\end{figure}
\vspace{-15pt}
\begin{figure}[H]
    \includegraphics[width=0.35\textwidth]{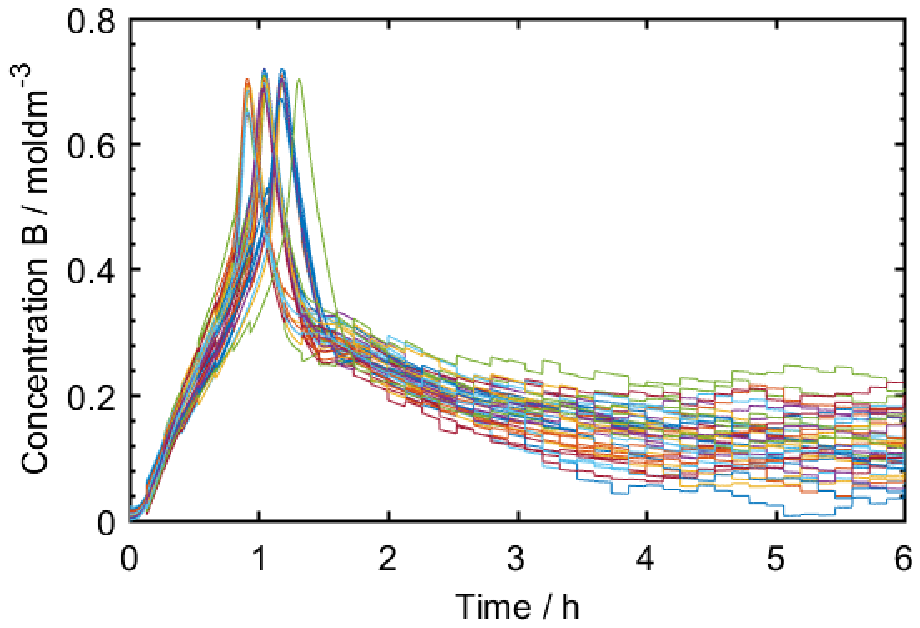}
    \centering
    \vspace{-10pt}
    \caption{Concentration B trajectories for 50 realizations}
\end{figure}
\vspace{-15pt}
\begin{figure}[H]
    \includegraphics[width=0.35\textwidth]{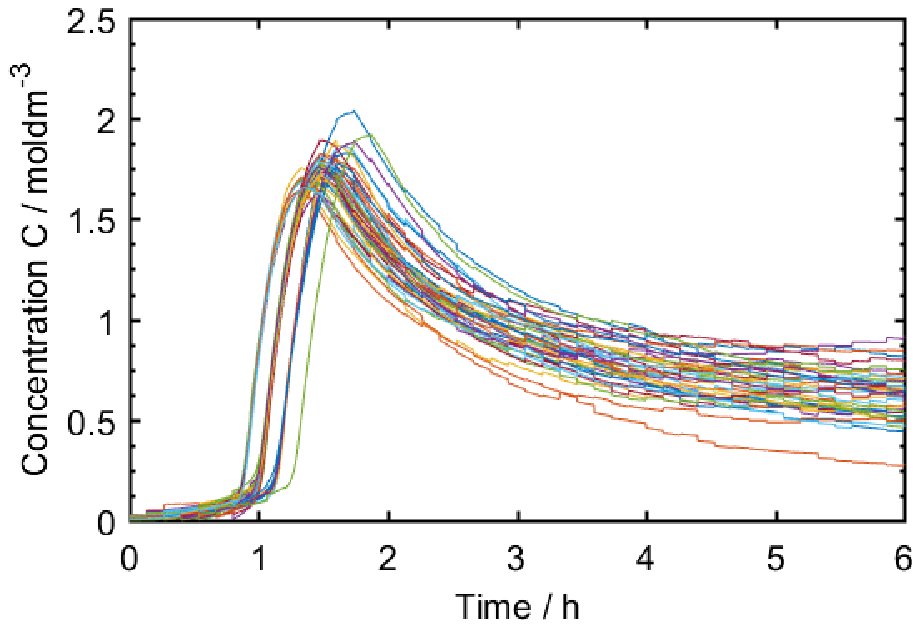}
    \centering
    \vspace{-10pt}
    \caption{Concentration C trajectories for 50 realizations}
\end{figure}
\vspace{-15pt}
\begin{figure}[H]
    \includegraphics[width=0.35\textwidth]{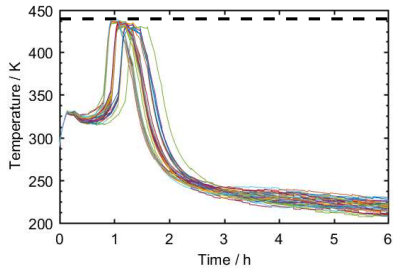}
    \centering
    \vspace{-10pt}
    \caption{Temperature trajectories for 50 realizations, path constraint shown as black dashed line}
\end{figure}
\vspace{-15pt}
\begin{figure}[H]
    \includegraphics[width=0.35\textwidth]{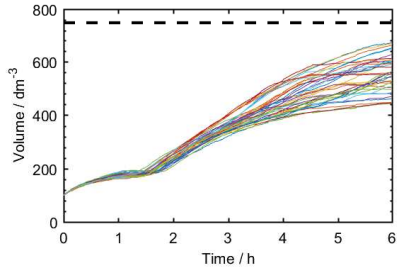}
    \centering
    \vspace{-10pt}
    \caption{Volume trajectories for 50 realizations, path constraint shown as black dashed line}
\end{figure}

\section{CONCLUSIONS}
In this paper, a SNMPC technique based on the incorporation of the UKF into the MPC problem was proposed. The approach uses the UKF to propagate both state estimate error and general additive uncertainty from disturbances forward in time. Linear joint chance constraints of the states could be easily implemented, since the states were assumed to follow a Gaussian distribution. The resulting OCP ensures feasibility through the probabilistic constraints. The concept of the "robust horizon" was introduced to handle the problem of growing covariances as the time horizon grows, which otherwise would lead to either a too conservative controller or infeasibilities in the OCP. The proposed UKF-SNMPC algorithm was then applied to a challenging semi-batch reactor case study with economic objective. Overall the algorithm was able to deal with the disturbances keeping the temperature of the reactor below the safety limit.        
 \\
\section{ACKNOWLEDGMENTS}
\begin{wrapfigure}{r}{0.16\textwidth}
\vspace{-5pt}
\includegraphics[width=0.16\textwidth]{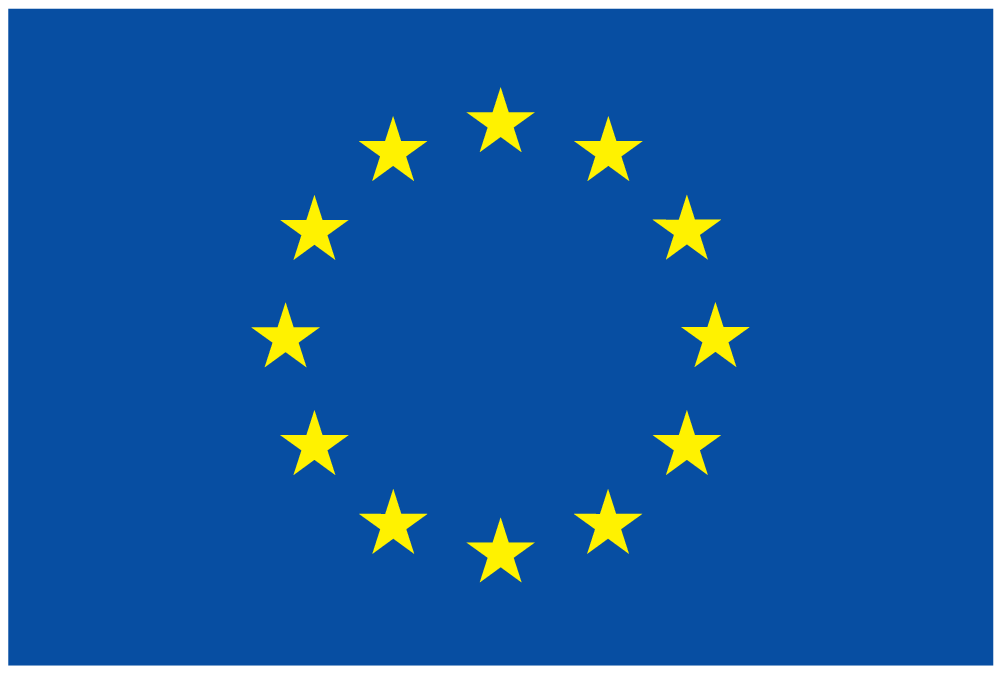}
\end{wrapfigure}
This project has received funding from the European Union's Horizon 2020 research and innovation programme under the Marie Sklodowska-Curie grant agreement No 675215.

\bibliographystyle{IEEEtran}
\bibliography{Stochastic_model_predictive_control}

\end{document}